\documentclass[10pt,a4paper,onecolumn]{article}
\usepackage{marginnote}
\usepackage{graphicx}
\usepackage{xcolor}
\usepackage{authblk,etoolbox}
\usepackage{titlesec}
\usepackage{calc}
\usepackage{tikz}
\usepackage{hyperref}
\hypersetup{colorlinks,breaklinks=true,
            urlcolor=[rgb]{0.0, 0.5, 1.0},
            linkcolor=[rgb]{0.0, 0.5, 1.0},
            citecolor = blue}
\usepackage{caption}
\usepackage{tcolorbox}
\usepackage{amssymb,amsmath}
\usepackage{ifxetex,ifluatex}
\usepackage{seqsplit}
\usepackage{xstring}
\usepackage{natbib}
\usepackage[T1]{fontenc}
\usepackage[utf8]{inputenc}
\usepackage{cleveref}
\usepackage{float}
\let\origfigure\figure
\let\endorigfigure\endfigure
\renewenvironment{figure}[1][2] {
    \expandafter\origfigure\expandafter[H]
} {
    \endorigfigure
}

\let\textttOrig=\texttt
\def\texttt#1{\expandafter\textttOrig{\seqsplit{#1}}}
\renewcommand{\seqinsert}{\ifmmode
  \allowbreak
  \else\penalty6000\hspace{0pt plus 0.02em}\fi}

\usepackage{abstract}

\makeatletter
\let\href@Orig=\href
\def\href@Urllike#1#2{\href@Orig{#1}{\begingroup
    \def\Url@String{#2}\Url@FormatString
    \endgroup}}
\def\href@Notdoi#1#2{\def\tempa{#1}\def\tempb{#2}%
  \ifx\tempa\tempb\relax\href@Urllike{#1}{#2}\else
  \href@Orig{#1}{#2}\fi}
\def\href#1#2{%
  \IfBeginWith{#1}{https://doi.org}%
  {\href@Urllike{#1}{#2}}{\href@Notdoi{#1}{#2}}}
\makeatother

\newlength{\cslhangindent}
\newlength{\csllabelwidth}
 {
  \setlength{\parindent}{0pt}
  \ifodd #1 \everypar{\setlength{\hangindent}{\cslhangindent}}\ignorespaces\fi
  \ifnum #2 > 0
  \setlength{\parskip}{#2\baselineskip}
  \fi
 }%
 {}
\usepackage{calc}

\usepackage[top=3.5cm, bottom=3cm, right=1.5cm, left=1.0cm,
            headheight=2.2cm, reversemp, includemp, marginparwidth=4.5cm]{geometry}

\titleformat{\section}
  {\normalfont\sffamily\Large\bfseries}
  {}{0pt}{}
\titleformat{\subsection}
  {\normalfont\sffamily\large\bfseries}
  {}{0pt}{}
\titleformat{\subsubsection}
  {\normalfont\sffamily\bfseries}
  {}{0pt}{}
\titleformat*{\paragraph}
  {\sffamily\normalsize}

\usepackage{fancyhdr}
\makeatletter
\let\ps@plain\ps@fancy
\fancyheadoffset[L]{4.5cm}
\fancyfootoffset[L]{4.5cm}

\definecolor{linky}{rgb}{0.0, 0.5, 1.0}

\newtcolorbox{repobox}
   {colback=red, colframe=red!75!black,
     boxrule=0.5pt, arc=2pt, left=6pt, right=6pt, top=3pt, bottom=3pt}

\newcommand{\ExternalLink}{%
   \tikz[x=1.2ex, y=1.2ex, baseline=-0.05ex]{%
       \begin{scope}[x=1ex, y=1ex]
           \clip (-0.1,-0.1)
               --++ (-0, 1.2)
               --++ (0.6, 0)
               --++ (0, -0.6)
               --++ (0.6, 0)
               --++ (0, -1);
           \path[draw,
               line width = 0.5,
               rounded corners=0.5]
               (0,0) rectangle (1,1);
       \end{scope}
       \path[draw, line width = 0.5] (0.5, 0.5)
           -- (1, 1);
       \path[draw, line width = 0.5] (0.6, 1)
           -- (1, 1) -- (1, 0.6);
       }
   }

\patchcmd{\@maketitle}{center}{flushleft}{}{}
\patchcmd{\@maketitle}{center}{flushleft}{}{}
\patchcmd{\@maketitle}{\LARGE}{\LARGE\sffamily}{}{}
\def\maketitle{{%
  
  \AB@maketitle}}
\makeatletter
\renewcommand\AB@affilsepx{ \protect\Affilfont}
\renewcommand\AB@affilnote[1]{{\bfseries #1}\hspace{3pt}}
\renewcommand{\affil}[2][]%
   {\newaffiltrue\let\AB@blk@and\AB@pand
      \if\relax#1\relax\def\AB@note{\AB@thenote}\else\def\AB@note{#1}%
        \setcounter{Maxaffil}{0}\fi
        \begingroup
        \let\href=\href@Orig
        \let\texttt=\textttOrig
        \let\protect\@unexpandable@protect
        \def\thanks{\protect\thanks}\def\footnote{\protect\footnote}%
        \@temptokena=\expandafter{\AB@authors}%
        {\def\\{\protect\\\protect\Affilfont}\xdef\AB@temp{#2}}%
         \xdef\AB@authors{\the\@temptokena\AB@las\AB@au@str
         \protect\\[\affilsep]\protect\Affilfont\AB@temp}%
         \gdef\AB@las{}\gdef\AB@au@str{}%
        {\def\\{, \ignorespaces}\xdef\AB@temp{#2}}%
        \@temptokena=\expandafter{\AB@affillist}%
        \xdef\AB@affillist{\the\@temptokena \AB@affilsep
          \AB@affilnote{\AB@note}\protect\Affilfont\AB@temp}%
      \endgroup
       \let\AB@affilsep\AB@affilsepx
}
\makeatother

\renewcommand\Affilfont{\sffamily\small\mdseries}
\ifnum 0\ifxetex 1\fi\ifluatex 1\fi=0 
  \usepackage[T1]{fontenc}
  \usepackage[utf8]{inputenc}

\else 
  \ifxetex
    \usepackage{mathspec}
    \usepackage{fontspec}

  \else
    \usepackage{fontspec}
  \fi
  \defaultfontfeatures{Ligatures=TeX,Scale=MatchLowercase}

\fi
\IfFileExists{upquote.sty}{\usepackage{upquote}}{}
\IfFileExists{microtype.sty}{%
\usepackage{microtype}
\UseMicrotypeSet[protrusion]{basicmath} 
}{}

\let\addcontentslineOrig=\addcontentsline
\def\addcontentsline#1#2#3{\bgroup
  \let\texttt=\textttOrig\addcontentslineOrig{#1}{#2}{#3}\egroup}
\let\markbothOrig\markboth
\def\markboth#1#2{\bgroup
  \let\texttt=\textttOrig\markbothOrig{#1}{#2}\egroup}
\let\markrightOrig\markright
\def\markright#1{\bgroup
  \let\texttt=\textttOrig\markrightOrig{#1}\egroup}

\usepackage{graphicx,grffile}
\makeatletter
\def\maxwidth{\ifdim\Gin@nat@width>\linewidth\linewidth\else\Gin@nat@width\fi}
\def\maxheight{\ifdim\Gin@nat@height>\textheight\textheight\else\Gin@nat@height\fi}
\makeatother
\setkeys{Gin}{width=\maxwidth,height=\maxheight,keepaspectratio}
\IfFileExists{parskip.sty}{%
\usepackage{parskip}
}{
\setlength{\parindent}{0pt}
\setlength{\parskip}{6pt plus 2pt minus 1pt}
}
\ifx\paragraph\undefined\else
\let\oldparagraph\paragraph
\renewcommand{\paragraph}[1]{\oldparagraph{#1}\mbox{}}
\fi
\ifx\subparagraph\undefined\else
\let\oldsubparagraph\subparagraph
\renewcommand{\subparagraph}[1]{\oldsubparagraph{#1}\mbox{}}
\fi

\DeclareFixedFont{\ttb}{T1}{txtt}{bx}{n}{10} 
\DeclareFixedFont{\ttm}{T1}{txtt}{m}{n}{10}  

\usepackage{color}
\definecolor{deepblue}{rgb}{0,0,0.5}
\definecolor{deepred}{rgb}{0.6,0,0}
\definecolor{deepgreen}{rgb}{0,0.5,0}

\usepackage{listings}

\newcommand\pythonstyle{\lstset{
language=Python,
basicstyle=\scriptsize,
morekeywords={self},              
keywordstyle=\scriptsize\color{deepblue},
emph={MyClass,__init__},          
emphstyle=\ttb\color{deepred},    
stringstyle=\color{deepgreen},
    numbers=left,
frame=tb,                         
showstringspaces=false
}}

\lstnewenvironment{python}[1][]
{
\pythonstyle
\lstset{#1}
}
{}

\newcommand\pythoninline[1]{{\pythonstyle\lstinline!#1!}}

\title{pylevin: efficient numerical integration of integrals containing up to three Bessel functions}

\date{\vspace{-7ex}}

\begin{document}
        \author[Argelander Institut für Astronomie]{Robert Reischke\thanks{reischke@posteo.net, rreischke@astro.uni-bonn.de}}
   
\maketitle

\marginpar{

  \begin{flushleft}
  \sffamily\small

  {\bfseries DOI:} \href{https://joss.theoj.org/papers/10.21105/joss.08618}{\color{linky}{10.21105/joss.08618}}

  \vspace{2mm}

  {\bfseries Software}
  \begin{itemize}
    \setlength\itemsep{0em}
    \item \href{N/A}{\color{linky}{Review}} \ExternalLink
    \item \href{https://github.com/rreischke/levin_bessel}{\color{linky}{Repository}} \ExternalLink
    \item \href{DOI unavailable}{\color{linky}{Archive}} \ExternalLink
  \end{itemize}

  \vspace{2mm}

  \par\noindent\hrulefill\par

  \vspace{2mm}

  {\bfseries Editor:} \href{https://nkrusch.github.io}{Neea Rusch} \ExternalLink \\
  \vspace{1mm}
    {\bfseries Reviewers:}
  \begin{itemize}
  \setlength\itemsep{0em}
    \item \href{https://github.com/mithun218}{@mithun218}
    \item \href{https://github.com/KumarSaurabh1992}{@KumarSaurabh1992}
    \end{itemize}
    \vspace{2mm}

  {\bfseries Submitted:} 24 February 2025\\
  {\bfseries Published:} 19 November 2025

  \vspace{2mm}
  {\bfseries License}\\
  Authors of papers retain copyright and release the work under a Creative Commons Attribution 4.0 International License (\href{http://creativecommons.org/licenses/by/4.0/}{\color{linky}{CC BY 4.0}}).

  \end{flushleft}
}

\vspace{.2cm}
 {
   {\footnotesize Argelander-Institut für Astronomie, Universität Bonn, Auf dem Hügel 71, D-53121 Bonn, Germany}}

\vspace{.5cm}
\begin{abstract}
Integrals involving highly oscillatory Bessel functions are notoriously challenging to compute using conventional integration techniques. While several methods are available, they predominantly cater to integrals with at most a single Bessel function, resulting in specialised yet highly optimised solutions. 
Here we present \texttt{pylevin}, a \texttt{python} package to efficiently compute integrals containing up to three Bessel functions of arbitrary order and arguments. The implementation makes use of Levin's method and allows for accurate and fast integration of these highly oscillatory integrals. In benchmarking \texttt{pylevin} against existing software for single Bessel function integrals, we find its speed comparable, usually within a factor of two, to specialised packages such as FFTLog. Furthermore, when dealing with integrals containing two or three Bessel functions, \texttt{pylevin} delivers performance up to four orders of magnitude faster than standard adaptive quadrature methods, while also exhibiting better stability for large Bessel function arguments. \texttt{pylevin} is available from source via
\href{https://github.com/rreischke/levin_bessel}{github} or directly from \href{https://pypi.org/project/pylevin/}{pypi}.
\end{abstract}

\section{Introduction}
\noindent
Bessel functions commonly arise in physical systems exhibiting rotational symmetry. Consequently, theoretical predictions for observables often require evaluating integrals over these functions, which are typically not analytically solvable and must be evaluated numerically. Standard integration methods, such as quadrature, are generally inefficient and unreliable when handling these integrals due to the rapid oscillations characteristic of Bessel functions. Thus, it is essential to develop general tools that can compute these challenging integrals efficiently and accurately. \texttt{pylevin} addresses this need by facilitating the calculation of the following frequently encountered types of integrals.

\begin{equation}
\label{eq:integral}
I_{\ell_1\ell_2\ell_3}(k_1,k_2,k_3) = \int_{a}^{b} \mathrm{d}x\,f(x) \prod_{i=1}^N \mathcal{J}_{\ell_i}(k_ix)\,,\quad N= 1,2,3\,,
\end{equation}
here $\mathcal{J}_\ell(x)$ denotes a spherical, $j_\ell$, or cylindrical Bessel function, $J_\ell$ of order $\ell$ and $f(x)$ can be any non-oscillatory function, i.e. that if the product, $\prod_{i=1}^n \mathcal{J}_{\ell_i}(k_ix)$, has a characteristic frequency $\omega_\mathcal{J}$ then the characteristic frequency of $f$ must satisfy $\omega_{f} \ll \omega_\mathcal{J}$ over the integration domain. Otherwise, the function $f$ can depend arbitrarily on the integration variable and the parameters $k_i$. 

In the literature, there exist different methods\footnote{The list of these excellent methods and packages here will be very biased and is driven by my own field of research.} to solve integrals of the type of \Cref{eq:integral}, with public implementations solving the integral up to $N = 2$. \citet{ogata_2005} proposed a method which has been, for example, implemented in the \texttt{hankel}\footnote{\href{https://hankel.readthedocs.io/en/latest/index.html}{https://hankel.readthedocs.io/en/latest/index.html}} package \citet{murray_2019} to calculate Hankel integrals and transforms. These are integrals with the parameters: $a=0$, $N=1$ and $\mathcal{J}_{\ell_1}(kx) = J_\nu (kx)$ being a cylindrical Bessel function of order $\nu$. The same problem is tackled using Fast Fourier Transforms (FFT) in logarithmic space \citep[FFTLog][]{talman_1978} implemented for example in \texttt{pyfftlog} package\footnote{\href{https://pyfftlog.readthedocs.io/en/latest/}{https://pyfftlog.readthedocs.io/en/latest/}} \citep{hamilton_2000}, implementing discrete versions of the Hankel transformations as well. Another FFTLog implementation can be found in the
\texttt{hankl} package\footnote{\href{https://github.com/minaskar/hankl}{https://github.com/minaskar/hankl}} \citep{karamanis_2021}. \citet{schonberg_2018}\footnote{\href{https://github.com/lesgourg/class_public}{https://github.com/lesgourg/class\_public}} and \citet{fang_2020}\footnote{\href{https://github.com/xfangcosmo/FFTLog-and-beyond}{https://github.com/xfangcosmo/FFTLog-and-beyond}} used FFTLog as well to compute angular power spectra in the non-Limber projection, which requires evaluating \Cref{eq:integral} over a compact range (which can formally extended to infinity) and for $N= 1$ and $\mathcal{J}_{\ell_1}(kx) = j_\ell(kx)$. Further applications of FFTLog were reported, for example, in \citet{fang_2017,assassi_2017,grasshorn_2018}, whereas \citet{campagne_2017,chiarenza_2024} rely on Chebyshev polynomials.

FFTLog was further generalised to two dimensions (2DFFTLog) to compute integrals over products of two Bessel functions in \citet{fang_2d-fftlog_2020}\footnote{\href{https://github.com/xfangcosmo/2DFFTLog}{https://github.com/xfangcosmo/2DFFTLog}}. A completely different method was suggested by \citet{levin_fast_1996,levin_analysis_1997,iserles_efficient_2005} (Levin's method) which relies on casting the integral into a solution of a differential equation. An adaptive version \citep[see e.g.][for details on convergence, especially for small $\omega_\mathcal{J}$]{chen_2022} of this method was used in a series of papers \citep{zieser_cross-correlation_2016,spurio_mancini_3d_2018,spurio_mancini_testing_2018} to calculate integrals up to $N=2$. 

All the implementations presented are tailored to very specific applications and are highly optimised. However, they may lack flexibility and can integrate at most two Bessel functions. This is exactly where \texttt{pylevin} enters the stage. It allows for solving the type of integrals presented in \Cref{eq:integral} and therefore includes all the integrals from the discussed methods and more. The code is written in C\texttt{++} but is wrapped in \texttt{python} using \texttt{pybind} and is based on implementations used in \cite{zieser_cross-correlation_2016} which is expanded and heavily optimised. We first discuss the implementation (for a description of the method, see the appendix) and then benchmark against some popular codes.

\section{Implementation}
\label{sec:implementation}
As in \citet{zieser_cross-correlation_2016}, we implement the algorithm in C\texttt{++} as an adaptive method, i.e. we calculate the estimate of \Cref{eq:integral} by constructing a solution over the whole integral between $[a,b]$ at using \Cref{eq:solution} with $n$ collocation points. This is followed by calculating a solution using $n/2$ collocation points. The interval $[a,b]$ is then recursively bisected in that sub-interval where the relative error, i.e. the relative difference between $n/2$ and $n$, is largest. We truncate the recursion once convergence is achieved. The key improvements of \texttt{pylevin} compared to this previous implementation are: 
\begin{enumerate}
    \item The inclusion of more integral types, flexibility and an easy high-level interface via \texttt{python}.
    \item Memory and speed optimisation as well as \texttt{openmp} parallelisation. 
    \item The possibility to reuse precomputed quantities. Specifically, the function $f(x)$ appears solely within the inhomogeneity of the linear system of equations (see Equation \ref{eq:solution}). As a result, the homogeneous solution can always be established for a specific type of integral and a given bisection. In many practical applications, it is common to compute a fixed type of integral from \Cref{eq:integral} multiple times, each time with slightly different $f(x)$. A typical example is running a Monte Carlo Markov Chain (MCMC) in cosmology (see the comparison with \texttt{hankl} and \text{pyCCL}), where only the matter power spectrum and radial weight function changes, but the general structure of the integrand is kept the same. Hence, \texttt{pylevin} can update the function $f(x)$ and reuse the solution to the homogeneous equation to construct the particular solution for the inhomogeneous solution from \Cref{eq:solution}. This speeds up calculations after the homogeneous solution is known by an order of magnitude.
\end{enumerate}
A basic example of running the code is provided in the code block below. For a comprehensive description, we refer to the \href{https://levin-bessel.readthedocs.io/en/latest/index.html}{readthedocs}\footnote{A \texttt{python} notebook with examples can be found \href{https://github.com/rreischke/levin_bessel/blob/main/tutorial/levin_tutorial.ipynb}{here} \ExternalLink}. Here, we will briefly outline the purpose of the various steps involved:

\begin{figure} 
\begin{python}
import pylevin as levin
import numpy as np

x = np.geomspace(1e-5,100,100)
number_integrands = 2
y = np.linspace(1,2, number_integrands) 
f_of_x = x[:,None]**(3*y[None,:]) + (x**2 +x)[:, None] 

integral_type = 0 
N_thread = 1
logx = True 
logy = True 
lp_single = levin.pylevin(integral_type, x, f_of_x, logx, logy, N_thread)

n_sub = 10 
n_bisec_max = 32
rel_acc = 1e-4 
boost_bessel = False 
verbose = False 
lp_single.set_levin(n_sub, n_bisec_max, rel_acc, boost_bessel, verbose)

k = np.geomspace(1e-3,1e4,1000)
ell = (5*np.ones_like(k)).astype(int)
a = x[0]*np.ones_like(k)
b = x[-1]*np.ones_like(k)
diagonal = False 
result_levin = np.zeros((len(k), number_integrands))
lp_single.levin_integrate_bessel_single(a, b, k, ell, diagonal, result_levin)

f_of_x_new = x[:,None]**(2.5*y[None,:]) + (x**1.5 +x)[:, None]
lp_single.update_integrand(x,f_of_x_new, logx, logy)
lp_single.levin_integrate_bessel_single(a, b, k, ell, diagonal, result_levin)
\end{python}
\end{figure}

\begin{itemize}
    \item \textbf{Line 4 - 6}: Defines the non-oscillatory part of the integrand, $f(x)$. The $x$-range should include the minimum $a$ and the maximum $b$ you want to use. Note that the code allows to pass many integrands at the same time, therefore \pythoninline{f_of_x} has the shape \pythoninline{(len(x), number_integrands)}. The second dimension always has to exist, even it is one. 
     \item \textbf{Line 9 - 13}: Defines the type of the integral considered (here a single spherical Bessel function), the number of \texttt{openmp} threads and how the interpolation of the integrand should be carried out. In the last line, an instance of the class is created. 
     \item \textbf{Line 15 - 20}: Defines the settings of the Levin collocation method. This is optional, as the class assumes default values for all parameters. Here we set the number of collocation points, the maximum number of bisections until the relative accuracy is reached, and whether the Bessel functions should be computed using the \texttt{GSL}\footnote{\href{https://www.gnu.org/software/gsl/}{https://www.gnu.org/software/gsl/}} or \texttt{BOOST}\footnote{\href{https://www.boost.org}{{https://www.boost.org}}}.
     \item \textbf{Line 22 - 28}: First we set the parameters, $a$, $b$, $\ell$ and $k$ at which the integral should be calculated. These must all be one-dimensional arrays of the same length. In this case, we evaluate the integral at $10^3$ different $k$ values for the same values of $a$, $b$ and $\ell$. If the number of integrands is equal to the number of parameters, one can choose to only calculate the diagonal of that square matrix. Here we compute all values and therefore allocate a \pythoninline{(len(k), number_integrand)} array to store the results. If \pythoninline{len(k) == number_integrand} and \pythoninline{diagonal = True}, the result will be an array of length \pythoninline{number_integrand}.
    \item \textbf{Line 30 - 32}: We can now define a new $f(x)$ over the same $x$-range and simply update the integrand. Calling the integrator once more utilises the previously constructed bisection, thereby significantly speeding up the calculation, as previously mentioned. However, there are a few caveats to consider: all parameters of the integral must remain unchanged. Furthermore, if the integrand alters sufficiently such that the previously employed bisection is no longer adequate for achieving convergence, the results may become erroneous. This issue may arise, for example, if the asymptotic behaviour of $f(x)$ varies dramatically. Nonetheless, for most practical applications, this is typically not a concern.
     
\end{itemize}

\section{Comparison with various codes}
\noindent
In this section, we compare \texttt{pylevin} with several other, more specialised, codes that mostly solve integrals over single Bessel functions. All reported computation times are averages over several runs. The comparison was done with an M3 processor with 8 cores (4 performance cores). For \texttt{pylevin}, we always measure time with the homogeneous solution of \Cref{eq:solution} known, as this is the most common case in practice and the most efficient way to run the code. 
\subsection{\texttt{hankel}}
\begin{figure}
    \centering
    \includegraphics[width=0.49\textwidth]{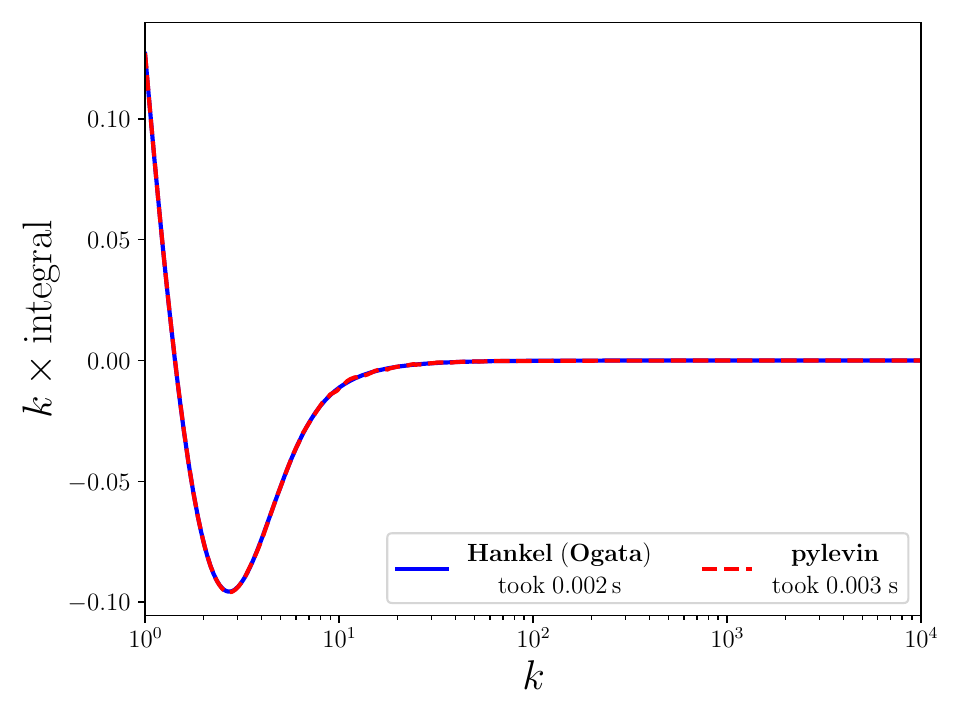}
    \includegraphics[width=0.49\textwidth]{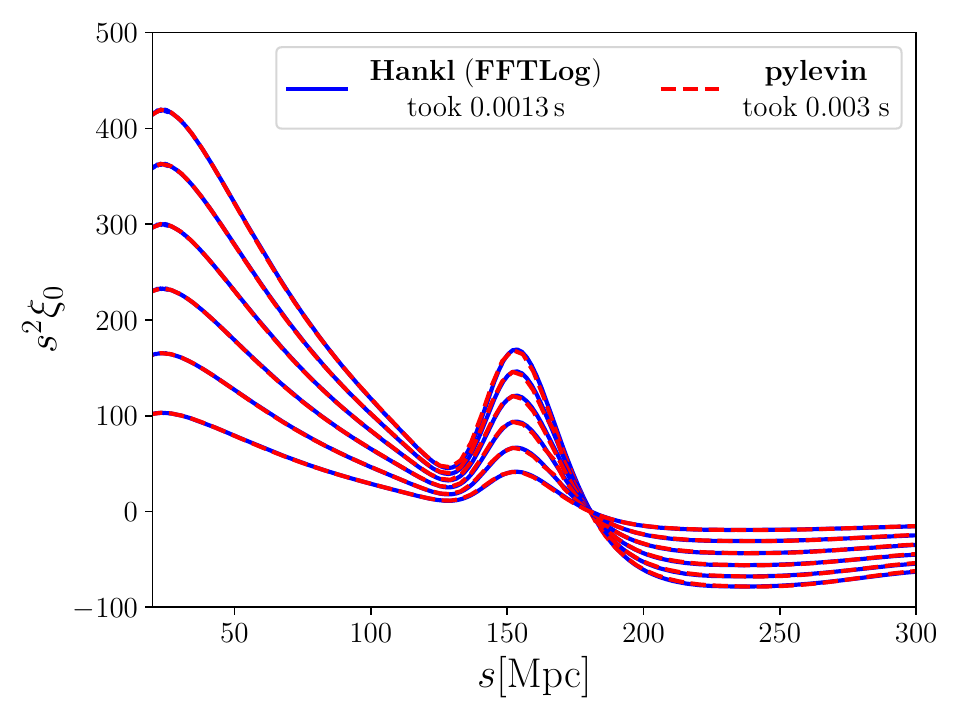}
    \caption{Comparison of \texttt{pylevin} with two methods to calculate a Hankel transformation. Dashed red is \texttt{pylevin} while solid blue is the alternative method. \textit{Left}: Integral in \Cref{eq:hankel} evaluated with the Ogata method using the \texttt{hankel} package. \textit{Right}: Integral in \Cref{eq:hankl} evaluated using the \texttt{hankl} package. Different lines refer to different redshifts.}
    \label{fig:hankel}
\end{figure}

\noindent
We calculate the following Hankel transformation:
\begin{equation}
\label{eq:hankel}
    \mathrm{integral}(k) = \int_0^\infty\frac{x^2}{x^2+1}J_0(kx)\;\mathrm{d}x\;,
\end{equation}
for 500 values of $k$ logarithmically-spaced between 1 and $10^4$. The result is depicted on the left side of \Cref{fig:hankel}. It can be seen that both methods agree very well and are roughly equally fast. While the Hankel transformation formally goes from zero to infinity, $a=10^{-5}$ and $b = 10^8$ were used for \texttt{pylevin}. This choice of course depends on the specific integrand.

\subsection{\texttt{hankl}}
\label{sec:hankl}
Here, we follow the cosmology example provided in the \texttt{hankl} documentation to compute the monopole of the galaxy power spectrum:
\begin{equation}
\label{eq:hankl}
    \xi_0(s) =\int_0^\infty\left(b^2+fb/3 +f^2/5\right)P_\mathrm{lin}(k)J_0(ks) k^2\;\mathrm{d}k\;,
\end{equation}
where $b$ is the galaxy bias, $f$ the logarithmic growth rate and $P_\mathrm{lin}(k)$ is the linear matter power spectrum, which is calculated using \texttt{camb} \citep{lewis_efficient_2000,lewis_cosmological_2002,howlett_cmb_2012} at six redshifts. Since \texttt{hankl} is FFT-based, it requires $k$ to be discretised; the FFT-dual will then be calculated at the inverse grid points. For this comparison, we use $2^{10}$ logarithmically-spaced points between $k=10^{-4}$ and $k = 1$ for the transformation to converge. For \texttt{pylevin}, the number of points where the transformation is evaluated is arbitrary. Here we use 100 points, which is more than enough to resolve all features in $\xi_0$.
The results are shown on the right of \Cref{fig:hankel} and good agreement can be found between the two methods, with \texttt{hankl} being roughly twice as fast as \texttt{pylevin}.

\subsection{\texttt{pyfftlog}}
\begin{figure}
    \centering
    \includegraphics[width=0.49\textwidth]{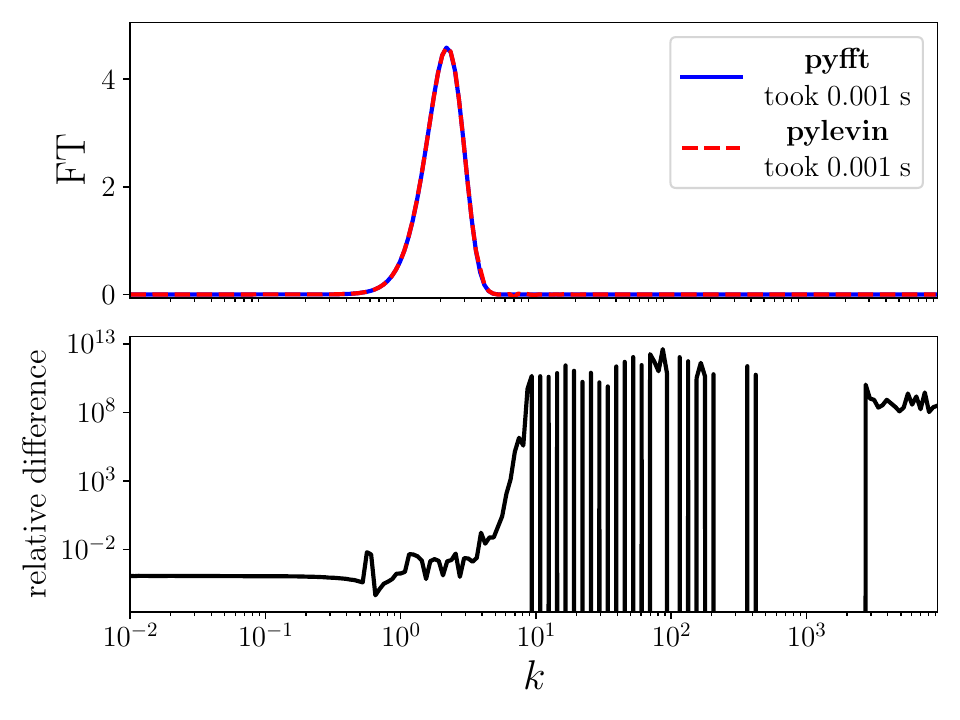}
    \includegraphics[width=0.49\textwidth]{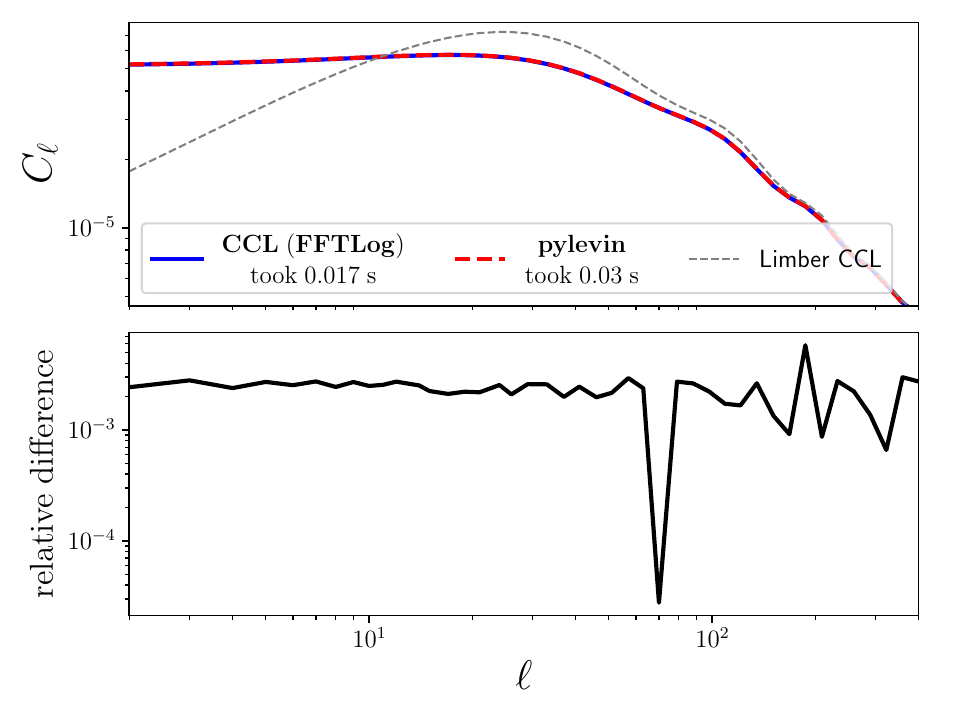}
    \caption{Comparison of \texttt{pylevin} with two other methods, the colour scheme is the same as in \Cref{fig:hankel}. \textit{Left}: Transformation defined in \Cref{eq:pyfftlog} with the \texttt{pyfftlog} package. We show the relative difference between the two methods in the lower panel. \textit{Right}: angular power spectra from \Cref{eq:nonlimber} computed with \texttt{pyCCL}.}
    \label{fig:ccl}
\end{figure}

For \texttt{pyfftlog} we use the following transformation:
\begin{equation}
\label{eq:pyfftlog}
    \mathrm{FT}(k) = \int_0^\infty r^{5}\mathrm{e}^{-r^2/2}J_4(kr)\;\mathrm{d}r\;.
\end{equation}
For \texttt{pyfftlog}, $2^8$ logarithmically-spaced points between $10^{-4}$ and $10^4$ for $r$ and hence also for $k$. \texttt{pylevin} is evaluated for the same number of points; this value could, however, be reduced due to the featureless transformation, thus increasing the speed. In the left panel of \Cref{fig:ccl}, the result of this exercise is shown. Good agreement between the two methods is found, with both taking the same amount of time. The large relative error at large values of $k$ is due to the small value of the integral, and hence purely numerical noise.

\subsection{\texttt{pyCCL}}
\label{sec:pyCCL}
Here, we compare the implementation of the non-Limber projection for the angular power spectrum:
\begin{equation}
\label{eq:nonlimber}
C_\ell = \frac{2}{\pi}\int\mathrm{d}{\chi_1} W(\chi_1)\int\mathrm{d}{\chi_2} W(\chi_2)\int k^2\mathrm{d}k\;  P_{\mathrm{m}}(k, \chi_1,\chi_2)
j_{\ell}(k\chi_1)j_{\ell}(k\chi_2) \;,
\end{equation}
for $W$, we assume a Gaussian shell in redshift with width $\sigma_z = 0.01$ centred at $z = 0.6$, $\chi(z)$ is the comoving distance. The matter power spectrum, $P_{\mathrm{m}}$, is again calculated with \texttt{camb}. The results from \texttt{pylevin} are compared to \texttt{CCL} \citep{chisari_core_2019}, which implements an FFTLog algorithm \citep{fang_2020, leonard_2023}. This implementation first solves the two integrals over $\chi_{1,2}$ using FFTLog and then carries out the remaining integral over $k$ and assumes that the $k$ and $\chi_{1,2}$ dependence in the matter power spectrum is separable. To be consistent, we follow the same approach and solve the $\chi_{1,2}$ integration using \texttt{pylevin} and then calculate the remaining $k$ integration with the composite Simpson's rule implemented in \texttt{scipy}. The right side of \Cref{fig:ccl} shows that the two methods agree very well with each other and that the method implemented in \texttt{CCL} is about a factor of 2 faster. If the power spectrum, however, were not separable on small $k$, as it can be the case in modified gravity scenarios, the \texttt{CCL} method would need to split the integral up into sub-intervals where the separability holds, slowing down the computation by a factor equal to the number of sub-intervals. This assumption is not done in \texttt{pylevin}. 

It should be noted that there is also an implementation in \texttt{julia}, which calculates \Cref{eq:nonlimber}, solving the $k$ integration first 
by expanding $P_\mathrm{m}$ into Chebyshev polynomials and precomputing the resulting integrals \citep{chiarenza_2024}. This was shown to be an order of magnitude faster after pre-computation, but it is also very tailored to this specific expression.

\section{Comparison with standard adaptive quadrature}
\begin{figure}
    \centering
    \includegraphics[width=0.49\textwidth]{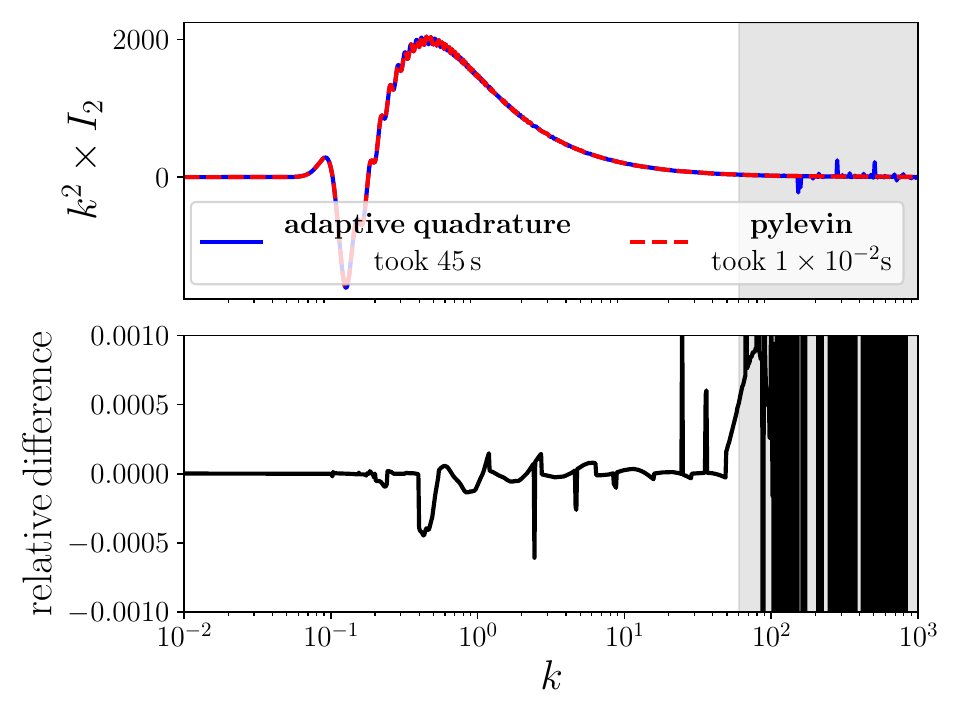}
        \includegraphics[width=0.49\textwidth]{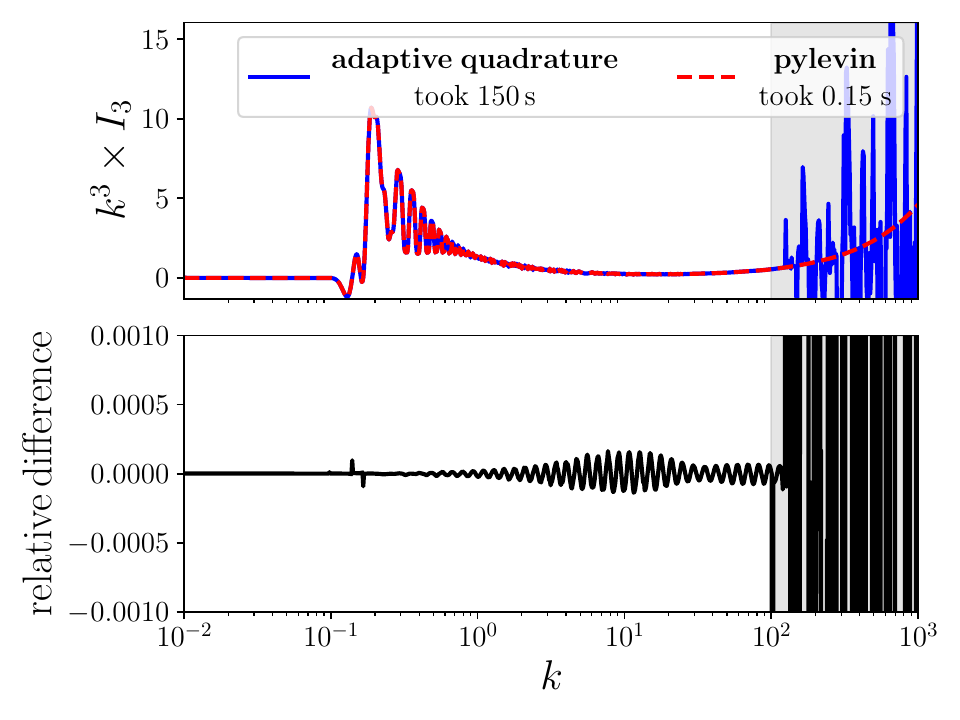}
    \caption{Comparison of \texttt{pylevin} with an adaptive quadrature (\texttt{scipy.integrate.quad}, the colour scheme is the same as in \Cref{fig:hankel}. \textit{Left}: Integral shown in \Cref{eq:two_bessel} with the relative difference in the lower panel. \textit{Right}: Integral shown in \Cref{eq:three_bessel} with the relative difference in the lower panel. The grey shaded area indicates the region where the adaptive quadrature fails to converge and catastrophically fails.}
    \label{fig:bm}
\end{figure}

Lastly, we calculate integrals over a product of two and three Bessel functions:
\begin{align}
\label{eq:two_bessel} 
I_2 = & \;\int_{10^{-5}}^{100} \mathrm{d}x \;(x^3 +x^2 +x)j_{10}(kx)j_5(kx)\;,
\\
\label{eq:three_bessel}
I_3 = & \;\int_{10^{-5}}^{100} \mathrm{d}x \;(x^3 +x^2 +x)j_{10}(kx)j_5(kx)j_{15}(kx)\;,
\end{align}
for $10^3$ $k$ values logarithmically spaced between $10^{-2}$ and $10^3$. In \Cref{fig:bm}, the comparison between the adaptive quadrature of \texttt{scipy} and $\texttt{pylevin}$ is shown with the integrals in \Cref{eq:two_bessel,eq:three_bessel} on the left and the right, respectively. In order for the quadrature to converge over an extended $k$-range, the number of maximum sub-intervals was increased to $10^3$ (from the default of 50). The grey-shaded area indicates where the quadrature fails to reach convergence even after this change. 
It is therefore clear that \texttt{pylevin} is more accurate and around three orders of magnitude faster than standard integration routines. 

\section{Summary}
We have introduced \texttt{pylevin}, a \texttt{python} implementation of Levin's method designed for the efficient evaluation of integrals involving up to three Bessel functions. The code was compared with methods that implement the special case of a single Bessel function using FFTLog or Ogata's method. Our tests show that \texttt{pylevin} is competitive with these more specialised codes and can reach comparable runtimes. For the case of integrals containing two and three Bessel functions, \texttt{pylevin} outperformed standard adaptive routines by almost four orders of magnitude.  \texttt{pylevin} is completely general with respect to the order and the arguments of the Bessel functions in the integral, and therefore presents an efficient and easy-to-use method to compute these challenging and time-consuming integrals.

\section{Acknowledgements}
I would like to thank Andrina Nicola and Paul Rogozenski for their help fixing issues with the \texttt{CCL} non-Limber comparison.
 
\bibliographystyle{mnras}
\bibliography{lib_paper} 

\newpage
\appendix

\section{Appendix: Levin's method}
\label{sec:levin}
\noindent
Levin's method \citep{levin_fast_1996,levin_analysis_1997} maps a quadrature problem of an oscillatory integral into the solution of a system of ordinary and linear differential equations. 
By assuming an orthonormal set of basis functions for the solution of the differential equation, the only thing which is required is to find the solution to a simple linear algebra problem. The method relies on integrals of the type:
\begin{align}\label{eq:i1}
I [\boldsymbol{f}(x,\{k_i\})] = \int_{x_1}^{x_2} \mathrm{d} x \,  \langle \boldsymbol{f} \left( x, \{k_i\} \right), {\bf{w}}(\{k_i\},x)\rangle \;.
\end{align}
where $\langle\cdot,\cdot\rangle$ denotes a scalar product. The non-oscillatory functions were bundled into a vector, $\boldsymbol{f}(x,\{k_i\})$ together with the oscillatory part denoted as $\boldsymbol{w}(x,\{k_i\})$. Both depend on the integration variable and a set of external parameters. For the method to work, the oscillatory part, $\boldsymbol{w}$, must satisfy:
\begin{equation}
    \frac{\mathrm{d}{\bf{w}}(x)}{\mathrm{d}x} = \boldsymbol{\sf{A}}(x) {\bf{w}}(x)\,,   
\end{equation}
with a matrix $\boldsymbol{\sf{A}}$. If one can find a vector $\boldsymbol{p}$ such that:
\begin{align}
\left\langle {\bf p}, {\bf w} \right\rangle' = \left\langle {\bf p}' + {\boldsymbol{\sf{A}}}^T{\bf p}, {\bf w} \right\rangle \approx \left\langle {\bf F}, {\bf w} \right\rangle\,,
\end{align}
the integral in Eq.~\ref{eq:i1} can be approximated by it.
Thus, the condition
$\left\langle {\bf p}, {\bf w} \right\rangle' = \left\langle {\bf F}, {\bf w}\right\rangle$, at $n$ collocation points $x_j , j = 1, 2, . . . , n$, leads to the following set of equations:
\begin{align}
  \left\langle {\bf p}' + {\boldsymbol{\sf{A}}}^T {\bf p} - {\bf F}, {\bf w} \right\rangle (x_j) = 0\,,
\end{align}
with the trivial solution:
\begin{align}\label{zeros}
{\bf p}'(x_j) + {\boldsymbol{\sf A}}^T (x_j) {\bf p}(x_j) = {\bf F}(x_j)\,.
\end{align}
Next, $\boldsymbol{p}$ is expanded into $n$ differentiable basis function $u_m$:
\begin{equation}
    {\bf p}(x) = \sum_{m=1}^{n} \boldsymbol{c}^{(m)}u_m(x)\;,
\end{equation}
where $\boldsymbol{c}^{(m)}$ are the coefficients of $\boldsymbol{p}$ in the function basis. Inserting the expansion back into the equation yields
\begin{equation}
\label{eq:solution}
    \sum_{m=1}^{n} \boldsymbol{c}^{(m)}u^\prime_m(x_j) + \sum_{m=1}^n \boldsymbol{\sf A}(x_j)\;\boldsymbol{c}^{(m)}u_m(x_j) = {\bf{F}}(x_j)\;,
\end{equation}
at the $n$ collocation points. For the basis functions, $u_m$, we choose the Chebyshev polynomials of the first kind, $T_n$ by first mapping the integral to the interval $[-1,1]$. This ensures more numerical stability than using simple polynomials, especially for low frequencies of the oscillatory part $\bf w$.

The explicit form of ${\bf{w}}$ and $\boldsymbol{\sf{A}}$ depends on the specific integral. In particular, we implement \Cref{eq:integral}, i.e. integrals with a single Bessel function, a product of two and a product of three Bessel functions. All types are implemented for spherical and cylindrical Bessel functions. We will now write down ${\bf{w}}$ and $\boldsymbol{\sf{A}}$ explicitly for the cases of spherical Bessel functions by using their recurrence relation.
These integrals can be evaluated by defining the following oscillatory parts
\begin{equation}
\begin{split}
   \bf{w}^{j}_\ell (x,k) = & \;\begin{pmatrix}
        j_{\ell}(xk) \\
        j_{\ell + 1}(xk)
        \end{pmatrix}\;,\\
    \bf{w}^{j}_{\ell_1,\ell_2} (x,k_1,k_2) = & \;\begin{pmatrix}
        j_{\ell_1}(xk_1) j_{\ell_2}(xk_2) \\
        j_{\ell_1 + 1}(xk_1) j_{\ell_2}(xk_2)\\
        j_{\ell_1}(xk_1)j_{\ell_2 + 1}(xk_2) \\
        j_{\ell_1 + 1}(xk_1)j_{\ell_2 + 1}(xk_2)
        \end{pmatrix} \;,   \\
        \bf{w}^{j}_{\ell_1,\ell_2,\ell_3} (x,k_1,k_2, k_3) = & \;\begin{pmatrix}
        j_{\ell_1}(xk_1) j_{\ell_2}(xk_2)j_{\ell_3}(xk_3) \\
        j_{\ell_1 + 1}(xk_1) j_{\ell_2}(xk_2)j_{\ell_3}(xk_3)\\
        j_{\ell_1}(xk_1)j_{\ell_2 + 1}(xk_2) j_{\ell_3}(xk_3)\\
         j_{\ell_1}(xk_1)j_{\ell_2}(xk_2) j_{\ell_3 + 1}(xk_3)\\
        j_{\ell_1 + 1}(xk_1)j_{\ell_2 + 1}(xk_2)j_{\ell_3}(k_3x)\\
        j_{\ell_1}(xk_1)j_{\ell_2 + 1}(xk_2)j_{\ell_3 + 1}(k_3x)\\
        j_{\ell_1 + 1}(xk_1)j_{\ell_2}(xk_2)j_{\ell_3 + 1}(k_3x)\\
        j_{\ell_1 + 1}(xk_1)j_{\ell_2 + 1}(xk_2)j_{\ell_3 + 1}(k_3x)
        \end{pmatrix}    \;.
\end{split}
\end{equation}
These would be exactly the same for cylindrical Bessel functions, $J_\ell$, by just replacing $j\to J$. For the matrix, $\boldsymbol{\sf{A}}^{j/J}$, one finds:
\begin{equation}
\begin{split}
         \boldsymbol{\sf{A}}^j_\ell (x,k)= &
    \begin{pmatrix}
        \frac{\ell}{x} & -k \\
        k & -\frac{\ell + 2}{x}
    \end{pmatrix}\;,\\
    \boldsymbol{\sf{A}}^J_\ell (x,k)= &
    \begin{pmatrix}
        \frac{\ell}{x} & -k \\
        k & -\frac{\ell + 1}{x}
    \end{pmatrix}\;,
\end{split}
\end{equation}
for $j$ and $J$ respectively. These can now be applied successively to obtain the matrices $\boldsymbol{\sf{A}}$ for the products of Bessel functions. From now on, we will only quote the results for $j$ since the ones for $J$ are easily obtained as well
\begin{equation}
   \begin{split}
       \boldsymbol{\sf{A}}^j_{\ell_1,\ell_2} (x,k_1,k_2)= &
    \begin{pmatrix}
        \frac{\ell_1 + \ell_2}{x} & -k_1 & -k_2 & 0  \\
        k_1 & \frac{\ell_2 - \ell_1  - 2}{x} & 0 &  - k_2 \\
        k_2 & 0 & \frac{\ell_1 - \ell_2  - 2}{x}& -k_1 \\
        0 &  k_2 & k_1 & -\frac{\ell_1 + \ell_2 +4}{x}
    \end{pmatrix}
    \end{split}\;,
\end{equation}
and
\begin{equation}
   \begin{split}
       \boldsymbol{\sf{A}}^j_{\ell_1,\ell_2,\ell_3} (x,k_1,k_2, k_3)= &
    \left(\begin{array}{cccc}
           \frac{\sum_i\ell_i}{x} & -k_1 & -k_2 & -k_3 \\
           k_1 & \frac{\ell_2  + \ell_3 - \ell_1  - 2}{x} & 0 & 0 \\
        k_2 & 0 & \frac{\ell_1 +\ell_3 -\ell_2 - 2}{x} & 0 \\
        k_3 & 0 & 0 &\frac{\ell_1 + \ell_2 - \ell_3 -2}{x} \\
                0 & k_2 & k_1 & 0 \\
        0 & k_2 & k_1 & 0 \\
        0 & 0 & k_3 & k_2 \\
        0 & k_3 & 0 & k_1  \\
        0 & 0 & 0 &0 \\
    \end{array}\right. \dots \\
     &\left.\begin{array}{cccc}
   0 & 0 & 0 & 0 \\
   -k_2 & 0 & -k_3 & 0\\
   -k_1 & -k_3 & 0 & 0 \\
    0 &-k_2 &-k_1 & 0 \\
    \frac{\ell_3 - \ell_1 - \ell_2 - 4}{x} & 0 & 0 & -k_3 \\
   0 & \frac{\ell_1 - \ell_2 - \ell_3 - 4}{x} &0 & -k_1\\
   0 & 0&  \frac{\ell_2 - \ell_1 -\ell_3 - 4}{x} & -k_2\\
   k_3 &k_1 & k_2 &-\frac{\sum_i\ell_i + 6}{x}
   \end{array}\right)
    \end{split}\;.
\end{equation}
In principle, it is thus easy to generalise the method to any other oscillatory function with some kind of recurrence relation. The user just needs to provide new functions for the function $\bf{w}$ and the matrix $\boldsymbol{\sf{A}}$ and define a new case in the code (that is, as long as the product of oscillatory functions does not extend beyond three).

\end{document}